\input amstex
\documentstyle{amsppt}
%----------------------------------------------------------------
% Title:     A note on the $\Sopfr(n)$ function.
% Author:    Ruslan Sharipov
% Comments:  AmSTeX, 7 pages, amsppt style
% MSC-class: 11N60, 11N64, 11-04
%----------------------------------------------------------------
%           Replacement for output macro definition
%
\catcode`@=11
\redefine\output@{%
  \def\break{\penalty-\@M}\let\par\endgraf
  \ifodd\pageno\global\hoffset=105pt\else\global\hoffset=8pt\fi  
  \shipout\vbox{%
    \ifplain@
      \let\makeheadline\relax \let\makefootline\relax
    \else
      \iffirstpage@ \global\firstpage@false
        \let\rightheadline\frheadline
        \let\leftheadline\flheadline
      \else
        \ifrunheads@ %\let\makefootline\relax
        \else \let\makeheadline\relax
        \fi
      \fi
    \fi
    \makeheadline \pagebody \makefootline}%
  \advancepageno \ifnum\outputpenalty>-\@MM\else\dosupereject\fi
}
\def\Beta{\mathchar"0\hexnumber@\rmfam 42}
\catcode`\@=\active
%----------------------------------------------------------------
\nopagenumbers
\def\negskp{\hskip -2pt}
\def\Sopfr{\operatorname{Sopfr}}
\def\blue#1{#1}

\catcode`#=11\def\diez{#}\catcode`#=6
\catcode`&=11\catcode`&=4
\catcode`_=11\catcode`_=8
%\catcode`~=11\def\volna{~}\catcode`~=\active
\def\mycite#1{\cite{\blue{#1}}\immediate\special{ps:
     ShrHPSdict begin /ShrBORDERthickness 0 def}}

\def\mytag#1{%
    \tag#1}
\def\mythetag#1{\thetag{\blue{#1}}\immediate\special{ps:
     ShrHPSdict begin /ShrBORDERthickness 0 def}}
\def\myrefno#1{\no#1}
\def\myhref#1#2{\blue{#2}\immediate\special{ps:
     ShrHPSdict begin /ShrBORDERthickness 0 def}}

\def\mytheorem#1{\csname proclaim\endcsname{Theorem #1}}
\def\mytheoremwithtitle#1#2{\csname proclaim\endcsname{Theorem #1#2}}

\def\mylemma#1{\csname proclaim\endcsname{Lemma #1}}
\def\mylemmawithtitle#1#2{\csname proclaim\endcsname{Lemma #1#2}}

\def\mycorollary#1{\csname proclaim\endcsname{Corollary #1}}

\def\myconjecture#1{\csname proclaim\endcsname{Conjecture #1}}
\def\myconjecturewithtitle#1#2{\csname proclaim\endcsname{Conjecture #1#2}}

%----------------------------------------------------------------
% Cyrillic fonts definition
%\font\eightcyr=wncyr8
%----------------------------------------------------------------
\pagewidth{360pt}
\pageheight{606pt}
\topmatter
\title
A note on the $\Sopfr(n)$ function.
\endtitle
\author
Ruslan Sharipov
\endauthor
\address Bashkir State University, 32 Zaki Validi street, 450074 Ufa, Russia
\endaddress
\email\myhref{mailto:r-sharipov\@mail.ru}{r-sharipov\@mail.ru}
\endemail
\abstract
    The $\Sopfr(n)$ function is defined as the sum of prime factors of $n$ 
each of which is taken with its multiplicity. This function is studied 
numerically. The analogy between $\Sopfr(n)$ and the primes distribution 
function is drawn and some conjectures for prime numbers formulated in terms 
of the $\Sopfr(n)$ function are suggested.
\endabstract
\endtopmatter
%\loadbold
%\loadeufb
\TagsOnRight
\document

\head
1. Introduction.
\endhead
     The $\Sopfr(n)$ function is defined as the sum of prime factors of its 
positive integer argument $n$ (see \mycite{1}). For $n=1$ this function is defined 
to be equal to zero: $\Sopfr(1)=0$. If $n$ is prime, then $\Sopfr(n)=n$. If $n$ 
is a product of prime numbers 
$$
\hskip -2em
n=p_1^{k_1}\cdot\ldots\cdot p_s^{k_s},
\mytag{1.1}
$$
then $\Sopfr(n)$ is calculated as the sum
$$
\hskip -2em
\Sopfr(n)=k_1\,p_1+\ldots+k_s\,p_s.
\mytag{1.2}
$$\par
     Note that the prime factors $p_1,\,\ldots,\,p_s$ in the sum \mythetag{1.2} are 
taken with their multiplicities $k_1,\,\ldots,\,k_s$ in the expansion \mythetag{1.1}.
Therefore the function $\Sopfr(n)$ is similar to the logarithm. One can easily
prove the following identity for it:
$$
\hskip -2em
\Sopfr(n_1\cdot n_2)=\Sopfr(n_1)+\Sopfr(n_2).
\mytag{1.3}
$$\par
     The $\Sopfr(n)$ function is used in defining Ruth-Aaron pairs named after 
two famous baseball players \myhref{http://en.wikipedia.org/wiki/Babe_Ruth}
{George Herman Ruth Jr\.} and \myhref{http://en.wikipedia.org/wiki/Hank_Aaron}
{Henry Louis Aaron} (see \mycite{2}). In mathematics a Ruth-Aaron pair is a pair 
of consecutive numbers $n$ and $n+1$ whose sums of prime factors are equal to
each other:
$$
\hskip -2em
\Sopfr(n)=\Sopfr(n+1).
\mytag{1.4}
$$
The numbers $714$ and $715$ constitute the most famous Ruth-Aaron pair.\par
     Let $x$ be an integer number and let $p$, $q$, $r$, and $s$ be four
numbers expressed through $x$ by the following four polynomials:
$$
\xalignat 2
&\hskip -2em
p=8\,x+5, &&q=48\,x^2+24\,x-1,\\
\vspace{-1.5ex}
&&&\mytag{1.5}\\
\vspace{-1.5ex}
&\hskip -2em
r=2\,x+1, &&s=48\,x^2+30\,x-1.
\endxalignat
$$
Using the formulas \mythetag{1.5}, one easily derives that 
$$
\xalignat 2
&\hskip -2em
p\,q+1=2^2\,r\,s,
&&p+q=2\cdot 2+r+s.
\mytag{1.6}
\endxalignat
$$
Due to \mythetag{1.6} and \mythetag{1.3}, if $p$, $q$, $r$, $s$ all are prime 
numbers, then the numbers $n=p\,q$ and $n+1=4\,r\,s$ constitute a Ruth-Aaron pair, 
i\.\,e\. they satisfy the equality \mythetag{1.4}. Schinzel's H-conjecture 
(see \mycite{3}, \mycite{4}, and \mycite{5}) implies that there are infinitely 
many integer numbers $x$ such that the numbers $p$, $q$, $r$, and $s$ given by the 
polynomials \mythetag{1.5} all are prime.\par
     In this paper we treat $\Sopfr(n)$ as an analog of the primes distribution
function $\pi(n)$. The value $\pi(n)$ of this function is defined as the number of 
positive primes less than or equal to $n$. Gauss and Legendre (see \mycite{6})
in 1792--1808 conjectured the following asymptotic behavior of the function 
$\pi(n)$:
$$
\hskip -2em
\pi(n)\sim \frac{n}{\ln(n)}\text{\ \ as \ }n\to\infty.
\mytag{1.7}
$$
In 1849 and in 1852 P.~L.~Chebyshev proved two propositions very close to 
\mythetag{1.7}. The proposition \mythetag{1.7} itself was proved in 1896 
by Hadamard \mycite{7} and Val\'ee Poussin \mycite{8}. See \mycite{9} for
the modern explanation of their proof.\par
     The main goal of this paper is to study the function $\Sopfr(n)$ numerically
and formulate some conjectures similar to \mythetag{1.7} for this function. 
\head
2. The averaged $\Sopfr(n)$ function.
\endhead
     The $\Sopfr(n)$ function is quite irregular. Looking at its graph (see
\mycite{1}), one can find that the values of $\Sopfr(n)$ resemble random
numbers. In order to make them more regular we average them over intervals
between two consecutive squares:
$$
\hskip -2em
A(n)=\sum^{(n+1)^2}_{i=n^2+1}\frac{\Sopfr(i)}{(n+1)^2-n^2}.
\mytag{2.1}
$$
The function $A(n)$ in \mythetag{2.1} is the averaged $\Sopfr(n)$ function. 
We study its values in two intervals $1\leqslant n\leqslant 998$ and 
$1000\leqslant n\leqslant 3161$. The graph of the function \mythetag{2.1} in 
the first interval $1\leqslant n\leqslant 998$ is shown in Fig\.~2.1. It is 
presented by a sequence of points whose coordinates are rendered in 
logarithmic scale, i\.\,e\. $A_n=(x_n,y_n)$, where $x_n=\ln(A(n))$ and 
$y_n=\ln(n)$.\par
     Looking at Fig\.~2.1 below, one can see that the points $A_n$ with 
$n\geqslant 122\approx e^{4.8}$ are approximated by a straight line. We write 
the equation of this straight line as 
$$
\hskip -2em
x=\alpha\,y+\beta.  
\mytag{2.2}
$$
In order to calculate the parameters $\alpha$ and $\beta$ in \mythetag{2.2} 
we use the root mean squares method. For this purpose we use the following
quadratic deviation function:
$$
\hskip -2em
F(\alpha,\beta)=\sum^{998}_{n=122}(x_n-\alpha\,y_n-\beta)^2. 
\mytag{2.3}
$$
The quadratic function \mythetag{2.3} \pagebreak has exactly one minimum point 
$\alpha=\alpha_{\sssize\text{min}}$, $\beta=\beta_{\sssize\text{min}}$. 
This minimum point is determined by the following linear equations:
$$
\xalignat 2
&\hskip -2em
\frac{\partial F(\alpha,\beta)}{\partial\alpha}=0,
&&\frac{\partial F(\alpha,\beta)}{\partial\beta}=0.
\mytag{2.4}
\endxalignat
$$
\vskip -5pt\hbox to 0pt{\kern -20pt \includegraphics{sopfr01.eps}\hss}
\vskip 260pt\noindent The function \mythetag{2.3} and the equations \mythetag{2.4} 
are computed numerically. 
\vadjust{\vtop to 265pt{\vskip 10pt\hbox to 0pt{\kern 
0pt \includegraphics{sopfr02.eps}\hss}\vfill\vphantom{111}}}Solving them, 
we find the numeric values of $\alpha$ and 
$\beta$ at the minimum of  the function $F(\alpha,\beta)$:
$$
\xalignat 2
&\hskip -2em
\alpha\approx 1.820,
&&\beta\approx -0.847.
\mytag{2.5}
\endxalignat
$$
Having calculated the constants \mythetag{2.5}, now we draw the graph of the 
deviation function $\delta(n)=\ln(A(n))-\alpha\,\ln(n)-\beta$ in logarithmic 
scale. The graph of the function $\delta(n)$ in Fig\.~2.2 is presented by a series
of points $A_n=(x_n,y_n)$, where $x_n=\ln(n)$ and $y_n=\delta(n)$. Looking at this
graph, we derive the following inequality for the deviation function $\delta(n)$:
$$
\hskip -2em
-\delta_1<\delta(n)<\delta_1\text{, \ where \ }\delta_1=0.15\text{ \ and \ }
122\leqslant n\leqslant 998.
\mytag{2.6}
$$\par
     The next interval is $1000\leqslant n\leqslant 3161$. The graph of the 
function \mythetag{2.1} in this interval is shown in Fig\.~2.3. Again it is 
presented by a sequence of points whose coordinates are rendered in logarithmic 
scale, i\.\,e\. $A_n=(x_n,y_n)$, where $x_n=\ln(A(n))$ and $y_n=\ln(n)$.
\vadjust{\vskip 5pt\hbox to 0pt{\kern 
-20pt \includegraphics{sopfr03.eps}\hss}\vskip 240pt}The graph 
in Fig\.~2.3 is also approximated by a straight line. This straight line is
given by the equation \mythetag{2.2}. The coefficients $\alpha$ and $\beta$ in
this case are calculated by solving the equations \mythetag{2.4} for the 
following quadratic deviation function, which is  similar to \mythetag{2.3}:
$$
\hskip -2em
F(\alpha,\beta)=\sum^{3161}_{n=1000}(x_n-\alpha\,y_n-\beta)^2. 
\mytag{2.7}
$$
The minimum of the function \mythetag{2.7} corresponds to the following 
values of $\alpha$ and $\beta$:
$$
\xalignat 2
&\hskip -2em
\alpha\approx 1.860,
&&\beta\approx -1.115.
\mytag{2.8}
\endxalignat
$$\par
     The sharpness of the approximation of $A(n)$ by the straight line in Fig\.~2.3
is expressed through the deviation function $\delta(n)=\ln(A(n))-\alpha\,\ln(n)-\beta$, 
\pagebreak where $\alpha$ and $\beta$ are given by the formulas \mythetag{2.8}:
$$
\hskip -2em
-\delta_2<\delta(n)<\delta_1\text{, \ where \ }\delta_2=0.1\text{ \ and \ }
1000\leqslant n\leqslant 3161.
\mytag{2.9}
$$
The inequalities \mythetag{2.9} are similar to the above inequalities \mythetag{2.6}.
They are derived by drawing the graph of the function $\delta(n)$. This graph is
shown in Fig\.~2.4 \vadjust{\vskip 5pt\hbox to 0pt{\kern 
0pt \includegraphics{sopfr04.eps}\hss}\vskip 290pt}below.\par
\head
3. Approximation conjectures. 
\endhead
     Note that the parameter $\alpha$ in \mythetag{2.8} is greater than $\alpha$ 
in \mythetag{2.5}. This mean that the slope of the straight line approximating
the graph of the function $A(n)$ slightly grows as $n\to\infty$. To take into
account this growth we replace the linear approximation in \mythetag{2.2} by a
nonlinear one. We choose the following formula for it
$$
\hskip -2em
x=\alpha\,y+\beta+\gamma\,\ln(y)+\lambda\,e^{-y}+\mu\,e^{-2\,y}.
\mytag{3.1}
$$
The choice of \mythetag{3.1} means that $A(n)$ is approximated by the formula
$$
\hskip -2em
A(n)\approx B\,n^\alpha\,(\ln n)^\gamma\,\exp\biggl(\frac{\lambda}{n}
+\frac{\mu}{n^2}\biggr)\text{, \ where \ }B=e^\beta.
\mytag{3.2}
$$
In order to find the optimal values of the parameters $\alpha$, $\beta$,
$\gamma$, $\lambda$, and $\mu$ for the approximation \mythetag{3.2} we apply
the root mean squares method. Instead of \mythetag{2.3} and \mythetag{2.7}
in this case we use the following deviation function:
$$
\hskip -2em
F=\sum_{n=4}^{3161}(x_n-\alpha\,y_n-\beta-\gamma\,\ln y_n
-\lambda\,e^{-y_n}-\mu\,e^{-2\,y_n})^2.
\mytag{3.3}
$$
Remember that $x_n=\ln(A(n)$ and $y_n=\ln(n)$ in \mythetag{3.3}. The optimal 
values of $\alpha$, $\beta$, $\gamma$, $\lambda$, and $\mu$ are determined
by solving the equations 
$$
\xalignat 5
&\frac{\partial F}{\partial\alpha}=0,
&&\frac{\partial F}{\partial\beta}=0,
&&\frac{\partial F}{\partial\gamma}=0,
&&\frac{\partial F}{\partial\lambda}=0,
&&\frac{\partial F}{\partial\mu}=0.
\qquad
\mytag{3.4}
\endxalignat
$$
The equations \mythetag{3.4} are similar to the equations \mythetag{2.4}. Here is
their solution:
$$
\xalignat 5
&\alpha\approx 2.001,
&&\beta\approx -0.047,
&&\gamma\approx -1.056,
&&\lambda\approx 1.187,
&&\mu\approx -2.240.
\qquad\quad
\mytag{3.5}
\endxalignat
$$\par
     The exponential factor with $\lambda$ and $\mu$ in \mythetag{3.2} is a
decreasing function of $n$. For this reason we consider the function
$$
B(n)=\frac{A(n)}{n^\alpha\,(\ln n)^\gamma}
$$
and draw its graph. Like the graph of $A(n)$, it is presented as a sequence of
\vadjust{\vskip 5pt\hbox to 0pt{\kern 
0pt \includegraphics{sopfr05.eps}\hss}\vskip 265pt}points:\par
     Looking at the graph in Fig\.~3.1, we can formulate the following
conjecture.
\myconjecturewithtitle{3.1}{ (weak $\Sopfr(n)$ conjecture)} There are four constants
$\alpha$, $\gamma$, $B_1$ and $B_2$ such that the averaged $\Sopfr$-function 
$A(n)$ in \mythetag{2.1} obey the inequalities
$$
B_1\,{n^\alpha\,(\ln n)^\gamma}\leqslant A(n)\leqslant 
B_2\,{n^\alpha\,(\ln n)^\gamma}\text{\ \ for all \ }n>1.
\mytag{3.6}
$$
\endproclaim
     The graph points in Fig\.~3.1 condense to a band as $n\to\infty$. Its width
is restricted by the constants $B_1$ and $B_2$ in \mythetag{3.6}. The width of this
band can vanish at infinity. For this option we can formulate the following
conjecture. 
\myconjecturewithtitle{3.2}{ (strong $\Sopfr(n)$ conjecture)} There are three 
constants $\alpha$, $\gamma$, and $B$ such that $B>0$ and the following condition 
is fulfilled:
$$
A(n)\sim B\,{n^\alpha\,(\ln n)^\gamma}\text{\ \ as \ }n\to\infty.
\mytag{3.7}
$$
\endproclaim
      Note that the constants $\alpha$ and $\gamma$ in \mythetag{3.5} are are very 
close to integer numbers. Therefore we can formulate another conjecture. 
\myconjecture{3.3} The constants $\alpha$ and $\gamma$ either in \mythetag{3.6}
or in \mythetag{3.7} are explicit numbers $\alpha=2$ and $\gamma=-1$.
\endproclaim
      The averaged $\Sopfr(n)$ function \mythetag{2.1} is similar to the
primes distribution function. The above formulas \mythetag{3.6} and \mythetag{3.7} 
are similar to the formula \mythetag{1.7}. 

\enddocument
\end